\documentclass[12pt, leqno]{amsart}
\usepackage{amscd,amsmath,amsopn,amssymb,amsthm}
\usepackage[english]{babel}
\usepackage{hyperref}
\DeclareMathOperator{\diag}{diag}
\DeclareMathOperator{\ad}{ad}
\DeclareMathOperator{\Ric}{Ric}

\DeclareMathOperator{\trace}{trace}

\DeclareMathOperator{\sgn}{sgn}

\renewenvironment{proof}[1][Proof]{\textbf{#1.} }
{\ \rule{0.5em}{0.5em}}

\newtheorem{theorem}{Theorem}
\newtheorem{pred}{Proposition}
\newtheorem{lemma}{Lemma}
\newtheorem{cor}{Corollary}

\newtheorem{remark}{Remark}

\newtheorem{conjecture}{Conjecture}
\renewcommand{\arraystretch}{1.4}

\begin{document}

\title[The signature of the Ricci curvature \dots on 4-dimensional Lie groups]
{The signature of the Ricci curvature\\ of  left invariant Riemannian metrics\\ on 4-dimensional Lie groups}

\author{A.G.~Kremlyov, Yu.G.~Nikonorov}

\address{Kremlyov Anton Gennadievich \newline
Rubtsovsk Industrial Institute
of Altai State \newline
Technical University after
I.I.~Polzunov
\newline  Rubtsovsk, Traktornaya, 2/6, 658207, Russia}

\email{kremant@mail.ru}

\address{Nikonorov Yurii  Gennadievich \newline
South Mathematical Institute
of  Vladikavkaz \newline
Scientific   Centre
of
the Russian Academy of Sciences,\newline
Vladikavkaz, Markus st. 22,
362027, Russia}

\email{nikonorov2006@mail.ru}

\begin{abstract}
In this paper, we present the classification of all possible signatures of the Ricci curvature of
left-invariant Riemannian metrics on 4-dimensional Lie groups and discuss some related questions.

\vspace{2mm} \noindent Key words and phrases: Riemannian
manifold, Lie groups, Ricci curvatures.
\end{abstract}

\maketitle

This brief note contains the main results of our papers \cite{KremNik} and \cite{KremNik1},
devoted to the classification of all possible signatures of the Ricci
curvatures of left-invariant metrics on 4-dimensional Lie groups.
We also discuss some related questions.

\section{The main results}

It is well known that various restrictions on the curvature of a Riemannian
manifold allow to get some information on its geometrical and topological properties.
Some results of this type are related to the Ricci curvature.
For example, the Myers theorem \cite{Mye} states that a complete Riemannian
manifold with positive Ricci curvature is compact and has finite fundamental group.
Note that the Ricci curvature is more informative for homogeneous Riemannian manifolds. For instance,
a Riemannian homogeneous manifold with negative Ricci curvature is noncompact by the
Bochner theorem \cite{Boch}.

For a homogeneous space $G/H$ (where $H$ is a compact subgroup of the Lie group $G$), it is naturally  to find
some general properties of the Ricci curvature of all $G$-invariant  Riemannian metrics on $G/H$.
One can make this problem more precise and concrete in various ways.
One variant is to consider the following question:
{\it What are possible signatures of the Ricci operators for all $G$-invariant Riemannian  metrics on the space $G/H$?}
There is a hope that this question could be completely resolved at least for the case
of low dimensions.
From the paper \cite{Mi} of J.~Milnor we know the answer to this question for dimensions $\leq3$.
The papers \cite{Ishihara, BB, Patr}
give us the corresponding answer for all four-dimensional homogeneous spaces, different from the Lie groups.

We recall briefly some structure results on homogeneous Riemannian
manifolds of low dimensions.
All $2$-dimensional simply connected homogeneous Riemannian manifolds are symmetric.
Each $3$-dimensional simply connected homogeneous Riemannian
manifold is either a symmetric space, or a Lie group supplied with a left-invariant Riemannian metric \cite{BB,Patr}.
As noted before, the signatures of the Ricci curvature of left invariant metrics on
$3$-dimensional Lie groups were classified in \cite{Mi}
(some refinements and generalizations one can find in \cite{KoNik, HaLee, RSch}).
The classification of $4$-dimensional homogeneous Riemannian
manifolds was obtained by S.~Ishihara \cite{Ishihara}.
G.~Jensen proved that every $4$-dimensional simply connected homogeneous Einstein manifold
is isometric to a symmetric space~\cite{Jen69}.
Some refinements of results of \cite{Ishihara} are obtained by L.~Berard-Bergery \cite{BB} and
V.~Patrangenaru \cite{Patr}.
In particular, it was shown in  \cite{Ishihara, BB, Patr} that an arbitrary simply connected homogeneous
Riemannian manifold of dimension
$\leq 4$ is either a symmetric space, or a Lie group with a suitable left invariant Riemannian metric.
Recall, that each symmetric space is a direct metric product of
Euclidian space and some irreducible symmetric spaces. Moreover, every irreducible symmetric space is
Einstein
(i.~e. it has constant Ricci curvature).
Therefore, the problem of determination of all possible signatures of the Ricci operators
of invariant metrics on a given symmetric space has an obvious solution.
The same problem for left invariant metrics on a given Lie group is more hard and interesting.

For dimensions $n \geq 5$ we have only partial results. It is necessary to mention
the paper \cite{DM} of I.~Dotti-Miatello, where Ricci signatures of left-invariant
Riemannian metrics on two-step solvable unimodular Lie groups are determined,
and the paper \cite{Krem} of A.G.~Kremlyov, where he determined all possible signatures of the Ricci operators
for left invariant metrics on nilpotent five-dimensional Lie groups.

The signature of a symmetric operator $A$, defined on a $n$-dimensional Euclidian space, is a $n$-tuple
($\sgn(\lambda_1),\, \sgn(\lambda_2), \dots , \sgn(\lambda_n)$), where $\lambda_1 \leq \dots \leq \lambda_n$
are the eigenvalues of the operator $A$ and $\sgn(x)$ means the sign of a (real) number $x$.
We list all possible signatures for the 4-dimensional case in Table 1.

\begin{table}[t]\label{sgn}
\begin{center} {\bf                                                    Table 1.}
\end{center}
{
\begin{center}
\begin{tabular}{|c|c||c|c||c|c|}
\hline \rule[0mm]{0.0 cm}{0 mm} {\normalsize{}} \rule[0 mm]{0.0
cm}{0 mm} &  \rule[0mm]{0.1 cm}{0 mm}\normalsize {Signature}
\rule[0 mm]{0.1 cm}{0mm} & \rule[0mm]{0.0 cm}{0 mm} {\normalsize{}} \rule[0 mm]{0.0
cm}{0 mm} &  \rule[0mm]{0.1 cm}{0 mm}\normalsize {Signature}
\rule[0 mm]{0.1 cm}{0mm} & \rule[0mm]{0.0 cm}{0 mm} {\normalsize{}} \rule[0 mm]{0.0
cm}{0 mm} &  \rule[0mm]{0.1 cm}{0 mm}\normalsize {Signature}
\rule[0 mm]{0.1 cm}{0mm}
\\
\hline\hline
$1$ & $(-,-,-,-)$ & $6$ & $(-,-,+,+)$ & $11$ & $(0,0,0,0)$\\
\hline
$2$ & $(-,-,-,0)$ & $7$ & $(-,0,0,0)$ & $12$ & $(0,0,0,+)$\\
\hline
$3$ & $(-,-,-,+)$ & $8$ & $(-,0,0,+)$ & $13$ & $(0,0,+,+)$\\
\hline
$4$ & $(-,-,0,0)$ & $9$ & $(-,0,+,+)$ & $14$ & $(0,+,+,+)$\\
\hline
$5$ & $(-,-,0,+)$ & $10$ & $(-,+,+,+)$ & $15$ & $(+,+,+,+)$\\
\hline
\end {tabular}
\end {center}
}
\end{table}

Any local problem related to the Ricci curvature of a left-invariant Riemannian metric $\rho$
on a (connected) Lie group $G$ can be reformulated in terms
of the corresponding metric Lie algebra. Recall,
that the left invariant metric $\rho$ defines an~inner product~$Q$
on the~Lie algebra~$\mathfrak{g}$ of~$G$ and vice versa: each
inner product~$Q$ on~$\mathfrak{g}$ induces a~left-invariant metric~$\rho$
on~$G$. There are useful formulas for the Ricci operator
of the metric Lie algebra $(\mathfrak{g}, Q)$ (see, e.~g. \cite{Al5} or \cite{Bes}).

\begin{table}[p]\label{liel}
\begin{center} {\bf                                     Table 2.}
\end{center}
{
%\large
%\scriptsize
\normalsize
\begin{center}
\begin{tabular}{|l|l|}
\hline \rule[0mm]{0.0 cm}{0 mm} {\normalsize{Lie algebra}} \rule[0 mm]{0.0
cm}{0 mm} &  \rule[0mm]{0.1 cm}{0 mm}\normalsize {Nonzero commutation relations}
\rule[0 mm]{0.1 cm}{0mm}
\\
\hline\hline
$4A_{1}$ & $$\\
\hline
$A_2\oplus2A_1$ & $[e_1,e_2]=e_2$\\
\hline
$2A_2$ & $[e_1,e_2]=e_2$, $[e_3,e_4]=e_4$\\
\hline
$A_{3,1}\oplus A_{1}$ & $[e_{2},e_{3}]=e_{1}$\\
\hline
$A_{3,2}\oplus A_{1}$ & $[e_1,e_3]=e_1$, $[e_2,e_3]=e_1+e_2$\\
\hline
$A_{3,3}\oplus A_{1}$ & $[e_1,e_3]=e_1$, $[e_2,e_3]=e_2$\\
\hline
$A_{3,4}\oplus A_{1}$ & $[e_{1},e_{3}]=e_{1}, [e_{2},e_{3}]=-e_{2}$\\
\hline
\begin{tabular}{l}
$A_{3,5}^{\alpha}\oplus A_{1}$,\\
$0<|\alpha|<1$
\end{tabular}
& $[e_1,e_3]=e_1$, $[e_2,e_3]=\alpha e_2$\\
\hline
$A_{3,6}\oplus A_{1}$ & $[e_{1},e_{3}]=-e_{2}, [e_{2},e_{3}]=e_{1}$\\
\hline
$A_{3,7}^{\alpha}\oplus A_{1}$,
$\alpha>0$
& $[e_1,e_3]=\alpha e_1-e_2$, $[e_2,e_3]=e_1+\alpha e_2$\\
\hline
$A_{3,8}\oplus A_{1}$ & $[e_{1},e_{2}]=-e_{3}, [e_{3},e_{1}]=e_{2}, [e_{2},e_{3}]=e_{1}$\\
\hline
$A_{3,9}\oplus A_{1}$ & $[e_{1},e_{2}]=e_{3}, [e_{3},e_{1}]=e_{2}, [e_{2},e_{3}]=e_{1}$\\
\hline
$A_{4,1}$ & $[e_{2},e_{4}]=e_{1}, [e_{3},e_{4}]=e_{2}$\\
\hline
$A_{4,2}^{\alpha}$, $\alpha\neq0$ & $[e_1,e_4]=\alpha e_1$, $[e_2,e_4]=e_2$, $[e_3,e_4]=e_2+e_3$\\
\hline
$A_{4,3}$ & $[e_1,e_4]=e_1$, $[e_3,e_4]=e_2$\\
\hline
$A_{4,4}$ & $[e_1,e_4]=e_1$, $[e_2,e_4]=e_1+e_2$, $[e_3,e_4]=e_2+e_3$\\
\hline
\begin{tabular}{l}
$A_{4,5}^{\alpha,\beta}$, $\alpha\beta\neq0$,\\
$-1\leq\alpha\leq\beta\leq1$
\end{tabular}
& $[e_1,e_4]=e_1$, $[e_2,e_4]=\alpha e_2$, $[e_3,e_4]=\beta e_3$ \\
\hline
$A_{4,6}^{\alpha,\beta}$, $\alpha\neq0$,
$\beta\geq0$
& $[e_1,e_4]=\alpha e_{1}$,
$[e_2,e_4]=\beta e_2-e_3$, $[e_3,e_4]=e_2+\beta e_3$ \\
\hline
$A_{4,7}$ & $[e_2,e_3]=e_1$, $[e_1,e_4]=2e_1$, $[e_2,e_4]=e_2$, $[e_3,e_4]=e_2+e_3$\\
\hline
$A_{4,8}$ & $[e_{2},e_{3}]=e_{1}, [e_{2},e_{4}]=e_{2}, [e_{3},e_{4}]=-e_{3}$\\
\hline
$A_{4,9}^{\beta}$,
$-1<\beta\leq1$
& $[e_2,e_3]=e_1$, $[e_1,e_4]=(1+\beta)e_1$, $[e_2,e_4]=e_2$, $[e_3,e_4]=\beta e_3$\\
\hline
$A_{4,10}$ & $[e_{2},e_{3}]=e_{1}, [e_{2},e_{4}]=-e_{3}, [e_{3},e_{4}]=e_{2}$\\
\hline
$A_{4,11}^{\alpha}$,
$\alpha>0$
&
\begin{tabular}{l}
$[e_2,e_3]=e_1$, $[e_1,e_4]=2\alpha e_1$, $[e_2,e_4]=\alpha e_2-e_3$,\\
$[e_3,e_4]=e_2+\alpha e_3$
\end{tabular}\\
\hline
$A_{4,12}$ & $[e_1,e_3]=e_1$, $[e_2,e_3]=e_2$, $[e_1,e_4]=-e_2$, $[e_2,e_4]=e_1$\\
\hline
\end {tabular}
\end {center}
}
\end{table}

Note that the set of inner products on a given $n$-dimensional Lie algebra  $\mathfrak{g}$
has dimension $n(n+1)/2$. This set could be reduced with using of the automorphism group of
$\mathfrak{g}$ (see e.g. \cite{Jen2}). In the $4$-dimensional case, after this reduction
we get the space of ``representative'' inner products of dimension $\leq 6=10-4$.
The classification, up to the action of the automorphism group,
of the inner products on all $4$-dimensional Lie algebras is obtained in the paper \cite{Fis}.
Therefore, for a given $4$-dimensional Lie group (Lie algebra),
our original problem could be reduced to the study of the signature
of some $4\times 4$ symmetric matrix (the matrix of the Ricci operator in a suitable orthonormal basis),
which depends on $\leq 6$ parameters.

In Table 2, we present the classification of $4$-dimensional Lie algebras obtained by G.M. Mubarakzyanov \cite{Mub}
(see also \cite{Patera, Fis, DG, ABDO}).
Now we can formulate our main results
(it should be noted that some partial results in this direction were obtained also in the paper \cite{Chen} of D.~Chen).

\begin{table}[t]\label{maint1}
\begin{center} {\bf   Table 3.}
\end{center}
{
%\large
%\scriptsize
\normalsize
\begin{center}
\begin{tabular}{|l|c|c|c|c|c|c|c|c|c|c|c|c|c|c|c|c|}
\hline \rule[0mm]{0.0 cm}{0 mm} {\normalsize{}} \rule[0 mm]{0.0
cm}{0 mm} &  \multicolumn{15}{|c|}{\normalsize{Signature}}\\
\hline\hline
\rule[0mm]{0.0 cm}{0 mm} {\normalsize{Lie algebra}}
\rule[0 mm]{0.0cm}{0 mm} &
$1$&$2$&$3$&$4$&$5$&$6$&$7$&$8$&$9$&$10$&$11$&$12$&$13$&$14$&$15$
\\ \hline\hline
$4A_{1}$ & $-$& $-$& $-$& $-$& $-$& $-$& $-$& $-$& $-$& $-$& $+$& $-$& $-$& $-$& $-$\\
\hline
$A_{3,1}\oplus A_{1}$ & $-$ & $-$& $-$& $-$& $+$& $-$& $-$& $-$& $-$& $-$& $-$& $-$& $-$& $-$& $-$\\
\hline
$A_{3,4}\oplus A_{1}$ & $-$ & $-$& $+$& $-$& $+$& $+$& $+$& $-$& $-$& $-$& $-$& $-$& $-$& $-$& $-$\\
\hline
$A_{3,6}\oplus A_{1}$ & $-$ & $-$& $+$& $-$& $+$& $+$& $-$& $-$& $-$& $-$& $+$& $-$& $-$& $-$& $-$\\
\hline
$A_{3,8}\oplus A_{1}$ & $-$ & $-$& $+$& $-$& $+$& $+$& $+$& $-$& $-$& $-$& $-$& $-$& $-$& $-$& $-$\\
\hline
$A_{3,9}\oplus A_{1}$ & $-$ & $-$& $+$& $-$& $+$& $+$& $-$& $+$& $+$& $+$& $-$& $+$& $-$& $+$& $-$\\
\hline
$A_{4,1}$ & $-$ & $-$& $+$& $-$& $+$& $+$& $-$& $-$& $-$& $-$& $-$& $-$& $-$& $-$& $-$\\
\hline
$A_{4,2}^{-2}$ & $-$ & $-$& $+$& $-$& $+$& $+$& $-$& $-$& $-$& $-$& $-$& $-$& $-$& $-$& $-$\\
\hline
\begin{tabular}{l}
$A_{4,5}^{\alpha,-1-\alpha}$, \\
$\alpha\in(-1,-\frac{1}{2})$
\end{tabular}
& $-$ & $-$& $+$& $-$& $+$& $+$& $+$& $-$& $-$& $-$& $-$& $-$& $-$& $-$& $-$\\
\hline
$A_{4,5}^{-1/2,-1/2}$ & $-$ & $-$& $-$& $-$& $+$& $-$& $+$& $-$& $-$& $-$& $-$& $-$& $-$& $-$& $-$\\
\hline
$A_{4,6}^{-2\beta,\beta}$ & $-$ & $-$& $+$& $-$& $+$& $+$& $+$& $-$& $-$& $-$& $-$& $-$& $-$& $-$& $-$\\
\hline
$A_{4,8}$ & $-$ & $-$& $+$& $-$& $+$& $+$& $-$& $-$& $-$& $-$& $-$& $-$& $-$& $-$& $-$\\
\hline
$A_{4,10}$ & $-$ & $-$& $+$& $-$& $+$& $+$& $-$& $-$& $-$& $-$& $-$& $-$& $-$& $-$& $-$\\
\hline
\end {tabular}
\end {center}
}
\end{table}

\begin{theorem}[\cite{KremNik}]\label{main}
Let $\mathfrak{g}$ be a unimodular $4$-dimensional Lie algebra from Table 2,\linebreak
$s$ be an arbitrary signature from Table 1. Then
$s$ can be realized as the signature
of the Ricci operator for some inner product on $\mathfrak{g}$ if and only if
there is the sign ``+'' in the entry of Table 3, corresponding to the Lie algebra
$\mathfrak{g}$ and the signature $s$.
\end{theorem}

\begin{cor}[\cite{KremNik}]\label{maim1}
The signatures
$(-,-,-,-)$, $(-,-,-,0)$, $(-,-,0,0)$, $(0,0,+,+)$ and $(+,+,+,+)$
are not the signatures of the Ricci operators
for left-invariant Riemannian metrics on $4$-dimensional univodular Lie groups.
\end{cor}

\begin{cor}[\cite{KremNik}]\label{main2}
Let $\mathfrak{g}$~be a $4$-dimensional unimodular Lie algebra, $Q$~be an arbitrary inner product
on $\mathfrak{g}$,
and $S$~be the scalar curvature of the metric Lie algebra $(\mathfrak{g}, Q)$.
Then the following assertions are true:

1) If $\mathfrak{g}$ commutative, then $S=0$;

2) If $\mathfrak{g}$ is isomorphic to $A_{3,6}\oplus A_{1}$, then $S=0$ or $S<0$;

3) If $\mathfrak{g}$ is isomorphic to $A_{3,9}\oplus A_{1}=su(2)\oplus \mathbb{R}$, then $S$ may have any sign;

4) For all other Lie algebras $\mathfrak{g}$ the inequality $S<0$ holds.
\end{cor}

\renewcommand{\arraystretch}{1.3}
\begin{table}[p]\label{maint2}
\begin{center} {\bf   Table 4.}
\end{center}
{%\large
%\scriptsize
\normalsize
\begin{center}
\begin{tabular}{|l|c|c|c|c|c|c|c|c|c|c|c|c|c|c|c|c|}
\hline \rule[0mm]{0.0 cm}{0 mm} {\normalsize{}} \rule[0 mm]{0.0
cm}{0 mm} &  \multicolumn{15}{|c|}{\normalsize{Signature}}\\
\hline\hline
\rule[0mm]{0.0 cm}{0 mm} {\normalsize{Lie algebra}}
\rule[0 mm]{0.0cm}{0 mm} &
$1$&$2$&$3$&$4$&$5$&$6$&$7$&$8$&$9$&$10$&$11$&$12$&$13$&$14$&$15$
\\ \hline\hline
$A_{2}\oplus 2A_{1}$ & $-$& $-$& $-$& $+$& $+$& $-$& $-$& $-$& $-$& $-$& $-$& $-$& $-$& $-$& $-$\\
\hline
$2A_{2}$ & $+$& $+$& $+$& $+$& $+$& $+$& $-$& $-$& $-$& $-$& $-$& $-$& $-$& $-$& $-$\\
\hline
$A_{3,2}\oplus A_{1}$ & $-$& $+$& $+$& $+$& $+$& $+$& $-$& $-$& $-$& $-$& $-$& $-$& $-$& $-$& $-$\\
\hline
$A_{3,3}\oplus A_{1}$ & $-$& $+$& $+$& $-$& $-$& $-$& $-$& $-$& $-$& $-$& $-$& $-$& $-$& $-$& $-$\\
\hline
$A_{3,5}^{\alpha}\oplus A_{1},\, \alpha\in(-1,0) $
& $-$& $-$& $+$& $-$& $+$& $+$& $-$& $-$& $-$& $-$& $-$& $-$& $-$& $-$& $-$\\
\hline
$A_{3,5}^{\alpha}\oplus A_{1},\, \alpha\in(0,1) $ & $-$& $+$& $+$& $+$& $+$& $+$& $-$& $-$& $-$& $-$& $-$& $-$& $-$& $-$& $-$\\
\hline
$A_{3,7}^{\alpha}\oplus A_{1}$ & $-$& $+$& $+$& $+$& $+$& $+$& $-$& $-$& $-$& $-$& $-$& $-$& $-$& $-$& $-$\\
\hline
$A_{4,2}^{\alpha},\, \alpha<0,\,\alpha\neq-2$ & $-$& $-$& $+$& $-$& $+$& $+$& $-$& $-$& $-$& $-$& $-$& $-$& $-$& $-$& $-$\\
\hline
$A_{4,2}^{\alpha},\, \alpha>0,\, \alpha\neq1$ & $+$& $+$& $+$& $+$& $+$& $+$& $-$& $-$& $-$& $-$& $-$& $-$& $-$& $-$& $-$\\
\hline
$A_{4,2}^{1}$ & $+$& $+$& $+$& $-$& $-$& $-$& $-$& $-$& $-$& $-$& $-$& $-$& $-$& $-$& $-$\\
\hline
$A_{4,3}$ & $-$& $-$& $+$& $-$& $+$& $+$& $-$& $-$& $-$& $-$& $-$& $-$& $-$& $-$& $-$\\
\hline
$A_{4,4}$ & $+$& $+$& $+$& $+$& $+$& $+$& $-$& $-$& $-$& $-$& $-$& $-$& $-$& $-$& $-$\\
\hline
$A_{4,5}^{\alpha,\alpha},\, \alpha\in[-1,-\frac{1}{2})$
& $-$& $-$& $+$& $-$& $-$& $-$& $-$& $-$& $-$& $-$& $-$& $-$& $-$& $-$& $-$\\
\hline
$A_{4,5}^{\alpha,\alpha},\,\alpha\in(-\frac{1}{2},0)$
& $-$& $-$& $-$& $-$& $-$& $+$& $-$& $-$& $-$& $-$& $-$& $-$& $-$& $-$& $-$\\
\hline
$A_{4,5}^{\alpha,1}$, $\alpha\in[-1,0)$ & $-$& $-$& $+$& $-$& $-$& $-$& $-$& $-$& $-$& $-$& $-$& $-$& $-$& $-$& $-$\\
\hline
$A_{4,5}^{\alpha,\alpha}$, $\alpha\in(0,1)$ & $+$& $+$& $+$& $-$& $-$& $-$& $-$& $-$& $-$& $-$& $-$& $-$& $-$& $-$& $-$\\
\hline
$A_{4,5}^{\alpha,1}$, $\alpha\in(0,1)$ & $+$& $+$& $+$& $-$& $-$& $-$& $-$& $-$& $-$& $-$& $-$& $-$& $-$& $-$& $-$\\
\hline
$A_{4,5}^{1,1}$ & $+$& $-$& $-$& $-$& $-$& $-$& $-$& $-$& $-$& $-$& $-$& $-$& $-$& $-$& $-$\\
\hline
\begin{tabular}{l}
$A_{4,5}^{\alpha,\beta}$, $\alpha\in[-1,0)$
\end{tabular}
& $-$& $-$& $+$& $-$& $+$& $+$& $-$& $-$& $-$& $-$& $-$& $-$& $-$& $-$& $-$\\
\hline
\begin{tabular}{l}
$A_{4,5}^{\alpha,\beta}$, $\alpha\in(0,1)$,\\ $\alpha\neq\beta$
\end{tabular}
& $+$& $+$& $+$& $+$& $+$& $+$& $-$& $-$& $-$& $-$& $-$& $-$& $-$& $-$& $-$\\
\hline
\begin{tabular}{l}
$A_{4,6}^{\alpha,\beta}$, $\alpha<0$,\\
$\beta>0$, $\alpha\neq-2\beta$
\end{tabular}
& $-$& $-$& $+$& $-$& $+$& $+$& $-$& $-$& $-$& $-$& $-$& $-$& $-$& $-$& $-$\\
\hline
\begin{tabular}{l}
$A_{4,6}^{\alpha,\beta}$, $\alpha>0$ , $\beta>0$
\end{tabular}
& $+$& $+$& $+$& $+$& $+$& $+$& $-$& $-$& $-$& $-$& $-$& $-$& $-$& $-$& $-$\\
\hline
$A_{4,6}^{\alpha,0}$ & $-$& $-$& $+$& $+$& $+$& $+$& $-$& $-$& $-$& $-$& $-$& $-$& $-$& $-$& $-$\\
\hline
$A_{4,7}$ & $+$& $+$& $+$& $+$& $+$& $+$& $-$& $-$& $-$& $-$& $-$& $-$& $-$& $-$& $-$\\
\hline
$A_{4,9}^{\beta},\, \beta\in(-1,-1/2)$ & $-$& $-$& $+$& $-$& $+$& $+$& $-$& $-$& $-$& $-$& $-$& $-$& $-$& $-$& $-$\\
\hline
$A_{4,9}^{-1/2}$ & $-$& $-$& $+$& $+$& $+$& $+$& $-$& $-$& $-$& $-$& $-$& $-$& $-$& $-$& $-$\\
\hline
$A_{4,9}^{\beta},\, \beta\in(-1/2,1)$ & $+$& $+$& $+$& $+$& $+$& $+$& $-$& $-$& $-$& $-$& $-$& $-$& $-$& $-$& $-$\\
\hline
$A_{4,9}^{1}$ & $+$& $+$& $+$& $-$& $-$& $-$& $-$& $-$& $-$& $-$& $-$& $-$& $-$& $-$& $-$\\
\hline
$A_{4,11}^{\alpha}$ & $+$& $+$& $+$& $+$& $+$& $+$& $-$& $-$& $-$& $-$& $-$& $-$& $-$& $-$& $-$\\
\hline
$A_{4,12}$ & $+$& $+$& $+$& $+$& $+$& $+$& $-$& $-$& $-$& $-$& $-$& $-$& $-$& $-$& $-$\\
\hline

\end {tabular}
\end {center}
}
\end{table}

\begin{theorem}[\cite{KremNik1}]\label{mainnonunim}
Let $\mathfrak{g}$ be a non-unimodular $4$-dimensional Lie algebra from Table~2,
$s$ be an arbitrary signature from Table 1. Then
$s$ can be realized as the signature
of the Ricci operator for some inner product on $\mathfrak{g}$ if and only if
there is the sign ``+'' in the entry of Table 4, corresponding to the Lie algebra
$\mathfrak{g}$ and the signature $s$.
\end{theorem}

\begin{remark}
In Table 4, the lines, corresponded to the Lie algebras $A_{4,9}^{\beta}$, $\beta\in(-1,1]$, are corrected,
because some results on these Lie algebras in \cite{KremNik1} are wrong
(see the next section for details).
\end{remark}

\begin{cor}[\cite{KremNik1}]\label{mainnonunim1}
The Ricci operator of any left-invariant Riemannian metric on every $4$-dimensional non-unimodular
Lie group has at least two negative eigenvalues.
For every signature from the list $(-,-,-,-)$, $(-,-,-,0)$, $(-,-,-,+)$, $(-,-,0,0)$, $(-,-,0,+)$ and $(-,-,+,+)$,
there are a $4$-dimensional non-unimodular Lie group $G$ and a left-invariant Riemannian metric
$\rho$ on $G$ such that the signature of the Ricci operator of $(G,\rho)$ coincides with the chosen one.
\end{cor}

The proofs and discussions of all above results can be found in the papers \cite{KremNik, KremNik1}.
Our main technical tools are some ideas and structure results from \cite{Al5, Mi, NikNik, N15}.
It should be noted that in order to study some signatures (for instance, $(-,-,0,0)$)
we have developed some special methods.

\medskip
Note also that in \cite{KremNik, KremNik1} we have found special (very useful for various computations) bases
for every inner product on every four-dimensional Lie algebra. These auxiliary results could be very helpful
in order to study other problems related to four-dimensional Lie groups with left-invariant Riemannian metrics.
This is confirmed e.~g. by papers \cite{GRS} and \cite{GS}.

\medskip

Note that the scalar curvature of any left-invariant Riemannian metric on a non-unimodular Lie group
is negative \cite{Jen2, Mi}. Thus the Ricci operator in this case
has at least one negative eigenvalue. Corollary \ref{mainnonunim} states that there are at least 2 such eigenvalues
in the $4$-dimensional case. The same is true for the dimensions $2$ and $3$ (see \cite{Mi}).
We supposed that it is a general property on non-unimodular solvable Lie groups
(for dimensions $\leq 4$ every non-unimodular Lie group is solvable):

\begin{conjecture}[\cite{KremNik1}]\label{twonullfour}
Let $\mathfrak{s}$ be a non-unimodular solvable Lie algebra of an arbitrary dimension.
Then for every inner product $Q$ on
$\mathfrak{s}$, the Ricci operator of the metric Lie algebra $(\mathfrak{s}, Q)$ has at least
two negative eigenvalues.
\end{conjecture}

This conjecture was confirmed for dimensions $\leq 4$ and
for non-unimodular metabelian Lie algebras in \cite{KremNik1},
for all non-unimodular solvable metric Lie algebras of dimension $\leq 6$ in \cite{Cheb1},
for all completely solvable Lie algebras in~\cite{NikCheb}, and
for all Lie algebras with six-dimensional two-step nilpotent derived algebras in~\cite{Abiev}.
The proof of this conjecture in full generality was obtain in the paper \cite{N16}.

\section{One example: the Lie algebras $A_{4,9}^{\beta}$}

In this section, we examine the signatures of the Ricci operator of all possible inner products on the Lie algebras $A_{4,9}^{\beta}$, $\beta\in (-1,1]$.
We have chosen this example because Proposition 22 in \cite{KremNik1} is incorrect. Below we give its correct version.

\begin{pred}

Only the signatures
$(-,-,-,+)$,  $(-,-,0,+)$, and $(-,-,+,+)$
(items 1, 5, and 6 from Table 1) are realized as signatures of the Ricci operators for all inner products
on the Lie algebra $A_{4,9}^{\beta}$ for any $\beta\in (-1,-1/2)$.

Only the signatures
$(-,-,-,+)$, $(-,-,0,0)$, $(-,-,0,+)$, and $(-,-,+,+)$
(items 3--6 from Table 1) are realized as signatures of the Ricci operators for all inner products
on the Lie algebra $A_{4,9}^{-1/2}$.

Only the signatures
$(-,-,-,-)$, $(-,-,-,0)$, $(-,-,-,+)$, $(-,-,0,0)$, $(-,-,0,+)$, and $(-,-,+,+)$
(items 1--6 from Table 1) are realized as signatures of the Ricci operators for all inner products
on the Lie algebra $A_{4,9}^{\beta}$ for any $\beta\in (-1/2,1)$.

Only the signatures
$(-,-,-,-)$, $(-,-,-,0)$, and $(-,-,-,+)$
(items 1--3 from Table 1) are realized as signatures of the Ricci operators for all inner products
on the Lie algebra $A_{4,9}^{1}$.
\end{pred}

The rest of this section consists of the proof of this proposition.
Let $A$ be a symmetric matrix, then by $A_{i_1,\dots,\,i_s}$ we denote its submatrix, which is obtained by deleting
all rows and columns with the numbers $i_1,\dots,\,i_s$.
\medskip

We have the following lemma for the Lie algebras $A_{4,9}^{\beta}$ (we suppose that $\beta \in (-1,1]$ is fixed).

\begin{lemma}[Lemma 18 in \cite{KremNik1}] \label{lemma19}
For any inner product $(\cdot,\cdot)$ on the Lie algebra Ëè $A_{4,9}^{\beta}$, there exists an
$(\cdot,\cdot)$-orthonormal basis $\{f_i\}$ with the following nonzero structure constants:
$$
\begin{array}{llll}
C_{1,4}^1=a(\beta+1),& C_{2,3}^1=b, & C_{2,4}^1=c, &
C_{2,4}^2=a,\\
$$
$$
C_{3,4}^1=d, & C_{3,4}^2=f(1-\beta), & C_{3,4}^3=a\beta\,,
\end{array}
$$
where $a, b, c, d, f \in \mathbb{R}$, $a > 0$ and $b > 0$.
Conversely, for all real numbers  $a>0$, $b>0$, $c$, $d$, and $f$, a Lie algebra with the above-mentioned (nonzero) structure constants
is isomorphic to the Lie algebra $A_{4,9}^{\beta}$.
\end{lemma}

Below we will study the Ricci operator $\Ric$ of the metric Lie algebra $(\mathfrak{g}=A_{4,9}^{\beta},(\cdot,\cdot))$, which is represented by  a symmetric
$(4\times4)$-matrix in the basis  $\{f_i\}$.
An explicit formula for $\Ric$ could be found e.~g. in \cite{Al5}:
\begin{equation*}
\Ric = -\frac{1}{2} \sum\limits_i \ad_{f_i}^{*} \ad_{f_i} + \frac{1}{4} \sum\limits_i \ad_{f_i}\ad_{f_i}^{*}-\frac{1}{2}B -(\ad_H)^{s},
\end{equation*}
where $A^*$ means the metric adjoint of an arbitrary operator $A$, $B$ is the Killing operator
(i.~e. $(BX,Y)=\trace(\ad(X)\cdot \ad(Y))$ for $X,Y \in \mathfrak{g}$),
the vector $H$ is determined by the equality $\trace \ad(X)= (X,H)$ for all $X \in \mathfrak{g}$,
$(\ad_H)^{s}=\frac12(\ad_H+\ad_H^*)$. Direct calculations  imply the following equality in our case:
{\tiny
$$
\Ric=\frac{1}{2}\left( {\begin{array}{rrrr}
b^2+c^2+d^2-4a^2(1+\beta)^2 & -ac\beta+df(1-\beta)-cl & -d(a+l) & 0 \\
-ac\beta+df(1-\beta)-cl & -2al-b^2-c^2+f^2(1-\beta)^2 & -cd-af(1-\beta)^2-fl(1-\beta) & bd \\
-d(a+l) & -cd-af(1-\beta)^2-fl(1-\beta) & -4a^2\beta(1+\beta)-b^2-d^2-f^2(1-\beta)^2 & -bc \\
0 & bd & -bc & -r
\end{array}}
 \right),
$$}
where $l=2a(1+\beta)>0$ and $r=4a^2(\beta^2+\beta+1)+c^2+d^2+f^2(1-\beta)^2 >0$.

It is clear that the forth diagonal elements of $\Ric$ is negative.
Now, consider the submatrices $\Ric_{1,2}$ and $\Ric_{1,3}$.
It is easy to check that
\begin{eqnarray*}
4\Bigl(\det(\Ric_{1,2})+\det(\Ric_{1,3})\Bigr)=
16(1+\beta)^2(\beta^2+\beta+1)a^4 \\
+ \Bigl(8(\beta ^{2} + \beta + 1)b^{2} + 4(3\beta + 2 +2\beta ^{2})(c^{2} +d^{2}) + 4f^{2}(\beta ^{2}-1)^{2}\Bigr)a^{2}  + c^{4}+ d^{4} \\
+ (2f^2(\beta-1)^2   + c^{2}+ d^{2} )b^{2} +
(f^2(\beta-1)^2 + 2d^{2} )c^{2}+ f^2(\beta-1)^2d^{2} >0.
\end{eqnarray*}
From this we get that at least one of the matrices $\Ric_{1,2}$ and $\Ric_{1,3}$ is negative definite. It implies that $\Ric$ has at least 2 negative eigenvalue.
Therefore, the signature of the Ricci operator of
the metric Lie algebra $(A_{4,9}^{\beta},(\cdot,\cdot))$ should be one of the signatures 1--6 in Table 1.

\begin{remark}
By Theorem 3 in \cite{KremNik1}, the Ricci operator of any four-dimensional nonunimodular metric Lie algebra
has at least two negative eigenvalues. See the end of the previous section for a discussion of this result and its generalizations.
\end{remark}

\subsection{The case $\beta=1$.}
In this special case, the matrix of the Ricci operator has the form
$$
\Ric = \frac{1}{2}\left( {\begin{array}{cccc}
b^2+c^2+d^2-16a^2 & -5ac & -5ad & 0 \\
-5ac & -8a^2-b^2-c^2 & -cd & bd \\
-5ad & -cd & -8a^2-b^2-d^2 & -bc \\
0 & bd & -bc & -12a^2-c^2-d^2
\end{array}}
 \right).
$$
The submatrix $\Ric_1$, that produced from $\Ric$ by deleting of the first row and the first column,
is negative definite. Therefore, the signatures 4--6  from Table 1 are impossible.
On the other hand, the signatures  1, 2, and 3 are the signatures of the Ricci operator for some
inner product, as shown in Table 5.
\begin{table}[h]
\begin{center} {\bf Table 5.}
\end{center}
\normalsize
\begin{center}
\begin{tabular}{|c|c|c|c|c|c|}
\hline \rule[0mm]{0.0 cm}{0 mm} {\normalsize{N of signature}} \rule[0
mm]{0.0 cm}{0 mm} &  \rule[0mm]{0.1 cm}{0 mm}\normalsize {
Signature}\rule[0 mm]{0.0 cm}{0 mm} &  \rule[0mm]{0.1 cm}{0
mm}\normalsize{$a$}\rule[0 mm]{0.0 cm}{0 mm} &  \rule[0mm]{0.1
cm}{0 mm}\normalsize{$b$}\rule[0 mm]{0.0 cm}{0 mm} &
\rule[0mm]{0.1 cm}{0 mm}\normalsize{$c$}\rule[0 mm]{0.0 cm}{0 mm}
&  \rule[0mm]{0.1 cm}{0 mm}\normalsize{$d$}\rule[0 mm]{0.0 cm}{0
mm}
\\
\hline
$1$ & $(-,-,-,-)$ & $1$ & $1$ & $0$ & $0$  \\
\hline
$2$ & $(-,-,-,0)$ & $1$ & $4$ & $0$ & $0$  \\
\hline
$3$ & $(-,-,-,+)$ & $1$ & $6$ & $0$ & $0$  \\
\hline
\end {tabular}
\end {center}
\end{table}

\subsection{The case $\beta\in (-1/2,1)$.}

Let  $p(t)$ be the characteristic polynomial of the matrix $2\Ric$ with $c=d=0$.
If in addition $a=1$ and $b=2(1+\beta)$, then
\begin{eqnarray*}
p(t)= t(t + 4(\beta ^{2} + \beta  + 1) + f^{2}(\beta-1)^{2})
\times \Bigl(t^2 +12(1+\beta)^2t
-(1-\beta)^{4}f^4\\
- (5\beta ^{2} + 6\beta  + 5)(1-\beta)^{2}f^2+16(2+\beta)(1+2\beta)(1+\beta)^{2}\Bigr).
\end{eqnarray*}
In particular, if $f=0$, then
$$
p(t)=t(t + 4(\beta ^{2} + \beta  + 1))
(t+ 4(1+\beta)(1+2\beta))(t+ 4(2+\beta)(1+\beta))\,,
$$
and we get the signature $(-,-,-,0)$ (recall that $\beta>-1/2$).
For any big enough $f$  we obviously get the signature $(-,-,0,+)$.
Hence, for some $f\in (0,\infty)$ we will have the signature $(-,-,0,0)$ too.

If $a=1$ and $b=3(1+\beta)$, then
\begin{eqnarray*}
p(t)=(t - 5(1+\beta)^{2})(t + 4(\beta ^{2} +\beta  + 1) + f^{2}(1-\beta)^2)
\times \Bigl(t^2 +22(1+\beta)^{2}t\\
- (1-\beta)^{4}f^{4} - (5\beta ^{2} + 6\beta  + 5)(1-\beta )^{2}f^{2} + (9\beta  + 13)(13\beta  + 9)(\beta  + 1)^{2}
\Bigr).
\end{eqnarray*}
It is easy get the signatures $(-,-,-,+)$ and
$(-,-,+,+)$ for suitable choice of $f$
(the first we get for $f=0$, the second for big enough $f$).

If $a=b=1$ and $f=0$, then
$$
p(t)=(t + (2\,\beta  + 3)\,(2\,\beta  + 1))(t + 4(\beta ^{2} + \beta  + 1))
(t+(2\,\beta  + 1)^{2})(t+5+4\beta)\,,
$$
and we get the signature $(-,-,-,-)$.

Hence, for $\beta\in (-1/2,1)$ the Ricci operator $\Ric$ may have any of the signatures 1--6 from the Table~1.

\subsection{The case $\beta\in (-1,-1/2]$.}
First, we prove the following

\begin{lemma}
If the Ricci operator $\Ric$ is non-positive defined, then $\beta=-1/2$ and the signature of $\Ric$ is $(-,-,0,0)$.
\end{lemma}

\begin{proof}
Consider the matrices
$$
D=\diag(- 3\beta , - 1 - 2\beta, 1 - \beta ,0), \quad
Q(t)=\diag \left(
1,
\left(\begin{array}{cc}
{\cos}(t) & {\sin}(t)\\
- {\sin}(t) & {\cos}(t) \\
\end{array}\right)
,1\right),
$$
and $D(t)=Q(t)DQ^{-1}(t)$. Since the eigenvalues $- 3\beta , - 1 - 2\beta, 1 - \beta ,0$ of $D(t)$  are non-negative,
$D(t)$ is non-negative defined. Therefore, $\trace\bigl(\Ric \cdot D(t)\bigr) \leq 0$ for all $t \in \mathbb{R}$.

Direct calculations show that $2\trace\bigl(\Ric \cdot D(t)\bigr) =h_1(t)+h_2(t)$, where
\begin{eqnarray*}
h_1(t)= (2 + \beta )\bigl({\cos}(t)c - {\sin}(t)d\bigr)^{2} - (1+2\beta)(c^{2} + d^{2}),\\
h_2(t)=4(1+\beta)\bigl((5 - \cos^{2}(t))\beta (1+\beta)  +{\cos}(2t)\bigr)a^{2}\\
- \sin(2t)(1 - \beta )(3+\beta)(2 + \beta )a f - (2 + \beta )(1 - \beta )^{2}{\cos}(2t) f^{2}.
\end{eqnarray*}
Note, that $h_1(t)$ and $h_2(t)$ don't depend on $b$.
Since $1+2\beta \leq 0$, then $h_1(t) \geq 0$ for all values of $t$.

Now, let us choose $t_0$ such that (recall that $a \neq 0$)
$$
\cos(2t_0)>0 \mbox{\,\,\,  and \,\,\,}
a\sin(2t_0)+f\cos(2t_0)=0\,.
$$
If we substitute $f=-a\sin(2t_0)/\cos(2t_0)$ in $h_2(t_0)$, then we get
$$
h_2(t_0)=\frac{4a^2(1+\beta)}{\cos(2t_0)}\Bigl( (9\cos^{2}(t_0) - 5)\beta(1+\beta) + 1 \Bigr).
$$
The minimal value of the function $\beta \mapsto (9\cos^{2}(t_0) - 5)\beta(1+\beta) + 1$ on $[-1,-1/2]$
may be achieved only at the points $\beta=-1$ and $\beta=-1/2$.
But in these points the function under consideration has the values $1$ and $\frac {9}{4}\sin^{2}(t_0)$ respectively.
Hence, its minimal value is non-negative and it is equal to $0$ if and only if $\sin(t_0)=0$,
$\beta=-1/2$ and $f=0$ simultaneously. If, in addition, $h_1(t_0)=0$, then $c=0$.

Therefore, for $\beta \in (-1,-1/2)$ we get $2\trace\bigl(\Ric \cdot D(t_0)\bigr) =h_1(t_0)+h_2(t_0)>0$, that impossible.
Hence for $\beta \in (-1,-1/2)$ the Ricci operator $\Ric$ has at least one positive eigenvalue.

If $\beta=-1/2$, $c=f=0$, then we have
$$
\Ric =\frac{1}{2}
\left(
{\begin{array}{cccc}
b^{2} + d^{2} - a^{2} & 0 &  - 2ad & 0 \\
0 &  - b^{2} - 2a^{2} & 0 & bd \\
 - 2ad & 0 &  - b^{2} - d^{2} + a^{2} & 0 \\
0 & bd & 0 &  - 3a^{2} - d^{2}
\end{array}}
\right)
$$
Consider two its submatrices
$$
\left(
{\begin{array}{cc}
b^{2} + d^{2} - a^{2} &  - 2ad  \\
 - 2ad & - b^{2} - d^{2} + a^{2}  \\
\end{array}}
\right) \mbox{\,\,\, and \mbox\,\,\,}
\left(
{\begin{array}{cc}
 - b^{2} - 2a^{2} &  bd \\
bd &  - 3a^{2} - d^{2}
\end{array}}
\right)\,.
$$
The first one has zero trace, hence, it has the signature $(0,0)$ for $b=a$ and $d=0$ or the
signature $(-,+)$ otherwise. The second submatrix is negative defined.
Therefore, $\Ric$ may have only signature $(-,-,-,+)$ or $(-,-,0,0)$ for $\beta=-1/2$, $c=f=0$.
In particular, $\Ric$ may not have the signatures $(-,-,-,-)$ or $(-,-,-,0)$ for $\beta=-1/2$.
\end{proof}
\medskip

We have proved that the signatures  $(-,-,-,+)$ and $(-,-,0,0)$  are admissible, but the signatures $(-,-,-,-)$ and $(-,-,-,0)$ are not  admissible
for $\Ric$ in the case $\beta=-1/2$ (see the end of the proof of Lemma 2). Now, it suffices to check the signatures $(-,-,+,+)$ and $(-,-,0,+)$.
Suppose that $c=d=0$ and $b=2a$, then
$$
\Ric =\frac{1}{2}\diag\left(3a^{2}, \left(
\begin{array}{cc}
 - 6\,a^{2} + {\displaystyle \frac {9}{4}} \,f^{2} &  -
{\frac {15}{4}}af  \\
 - {\frac {15}{4}} af &  - 3a^{2} -
{\frac {9}{4}}f^{2} \\
\end{array} \right), - 3a^{2} - \frac {9}{4} f^{2}\right)\,.
$$
The trace of the depicted $(2\times2)$-submatrix is negative. It is easy to prove that
the signatures $(-,-,+,+)$ and $(-,-,0,+)$ are realized for some suitable $a$ and $f$.

Finally, we get that for $\beta=-1/2$ the Ricci operator $\Ric$ may have exactly one of the signatures
$(-,-,-,+)$, $(-,-,0,0)$, $(-,-,0,+)$, and $(-,-,+,+)$.

\medskip

Now we finish the study of the case $\beta\in (-1,-1/2)$. We have proved (see Lemma~2)
that the signatures $(-,-,-,-)$, $(-,-,-,0)$, and $(-,-,0,0)$ are impossible.
Now we check the signature $(-,-,-,+)$, $(-,-,0,+)$, and $(-,-,+,+)$.
Put $c=d=0$ and $b=2a$, then
$$
\Ric=\frac{1}{2}
\diag  \Bigl(- 4a^{2}\beta(2 + \beta ),\, A\,, - 4a^{2}(1+\beta+\beta^2) - f^{2}(1-\beta)^2
 \Bigr),
$$
where
$$
A=\left(
\begin{array}{cc}
 - 4a^{2}(2+\beta) + f^{2}(1-\beta)^2
 & -af(3+\beta)(1-\beta)
 \\
-af(3+\beta)(1-\beta) &  - 4a^{2}(1+\beta+\beta^2) - f^{2}(1-\beta)^2
\end{array}\right).
$$
It is easy to see that $\trace(A)=- 4a^{2}(3+2\beta+\beta^2)<0$ and
$$
\det(A)=-(1-\beta)^{4}f^{4} - a^{2}(5\beta ^{2}+ 6\beta  + 5)(1-\beta)^{2}f^{2} + 16a^{4}(2 +\beta )(1 + \beta  + \beta ^{2}).
$$
It is clear that $\det(A)$ can have any sign for suitable $a$ and $f$.
Therefore, for $\beta\in (-1,-1/2)$ the Ricci operator $\Ric$ may have exactly one of the signatures
$(-,-,-,+)$, $(-,-,0,+)$, and $(-,-,+,+)$.

\bigskip

The authors are indebted to  Nurlan~Abiev and  Yuri~Nikolayevsky
for helpful discussions concerning this paper.

\bigskip

\bibliographystyle{amsunsrt}

\bigskip

\end{document}